\documentclass[11pt]{article}
\usepackage{amsmath}
\usepackage{amssymb}
\usepackage{amsfonts}
\usepackage{eucal}
\usepackage{graphicx,url}
\usepackage{amssymb,amsmath, amsthm}
\usepackage{epstopdf}
\newtheorem{Theorem}{Theorem}[section]
\newtheorem{Proposition}[Theorem]{Proposition}

\newtheorem{Lemma}[Theorem]{Lemma}
\newtheorem{Remark}[Theorem]{Remark}
\begin{document}
\title{Asymptotic results for random flights}
\author{Alessandro De Gregorio\thanks{Dipartimento di Scienze
Statistiche, Sapienza Universit\`{a} di Roma, Piazzale Aldo Moro
5, I-00185 Rome, Italy. e-mail:
\texttt{alessandro.degregorio@uniroma1.it}}\and Claudio
Macci\thanks{Dipartimento di Matematica, Universit\`a di Roma Tor
Vergata, Via della Ricerca Scientifica, I-00133 Rome, Italy.
e-mail: \texttt{macci@mat.uniroma2.it}}}
\date{}
\maketitle
\begin{abstract}
\noindent The random flights are (continuous time) random walks
with finite velocity. Often, these models describe the stochastic
motions arising in biology. In this paper we study the large time
asymptotic behavior of random flights. We prove the large
deviation principle for conditional laws given the number of the
changes of direction, and for the non-conditional laws of some
standard random flights.
\end{abstract}

\section{Introduction}\label{sec:introduction}
The random flight models have been introduced in order to describe
the real motions and they can be connected with the biological
locomotion. The original formulation of the random flight problem
is due to Pearson, who considers a random walk with fixed and
constant steps. The Pearson's model deals with a random walker
moving in the plane in straight lines with fixed length and
turning through any angle whatever. This random walk has been
introduced by Pearson for modeling the migration of mosquitoes
invading cleared jungles. Over the years the random flights have
attracted the attention of different researchers and the original
model has been extended randomizing the displacements of the walk;
see for example \cite{S}, \cite{OD}, \cite{LeC}, \cite{PRD},
\cite{DO}.

Several experiments have highlighted that many aspects of the cell
biology are governed by stochastic motions. For instance, some
motile bacteria, such as {\it Escherichia coli} and {\it
Salmonella typhimurium}, and blood cells have motion which is
well-approximated by random straight-lines with discrete changes
of direction (see \cite{NO}, and references therein). Therefore,
the random flights seem to be suitable for the analytic
characterization of the behavior of these micro-organisms.

\begin{figure}[ht]
\begin{center}
\includegraphics[angle=0,width=0.5\textwidth]{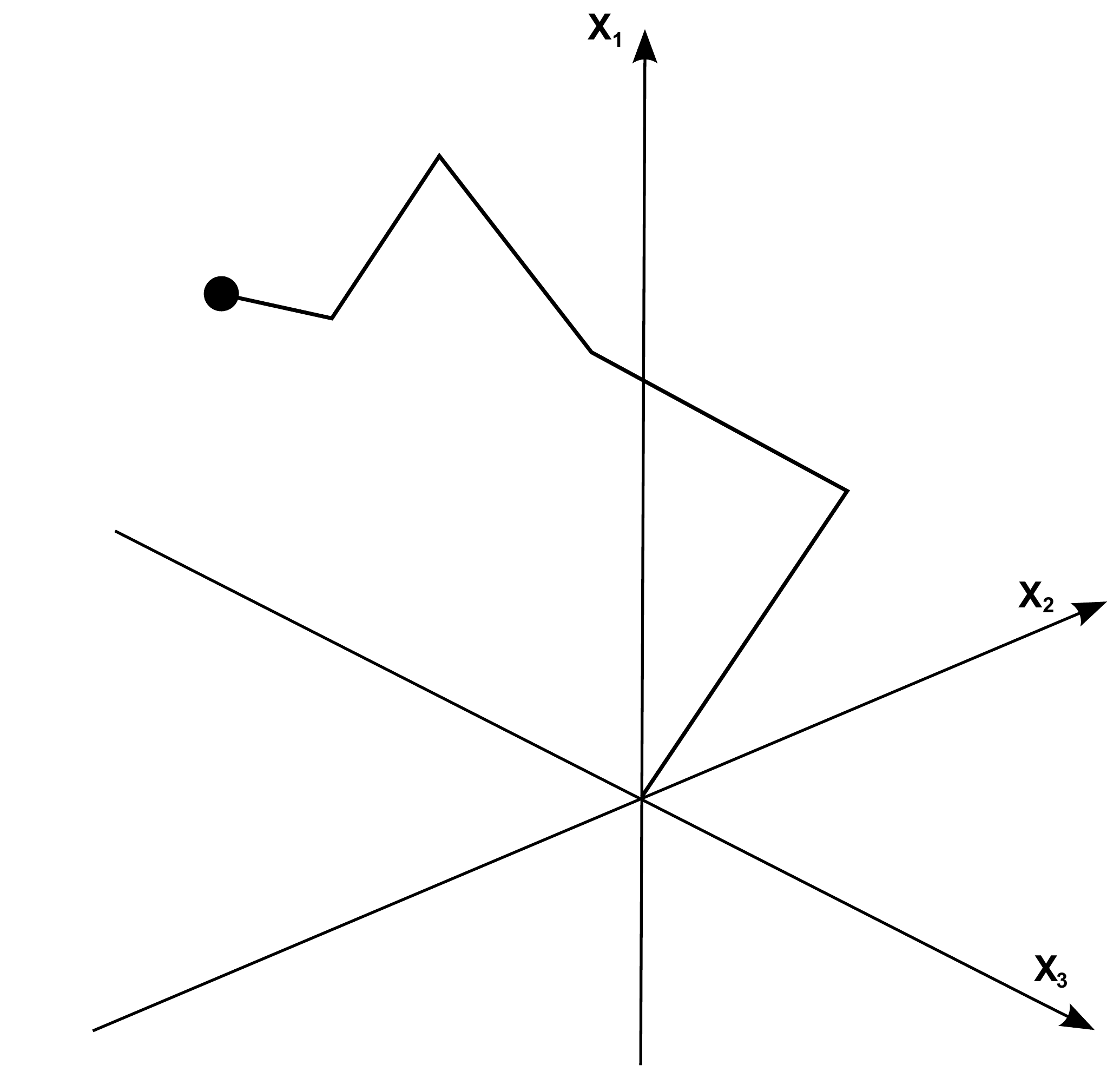}
\caption{A sample path in $\mathbb{R}^3$ consisting of $n=4$
changes of direction.}\label{fig1}
\end{center}
\end{figure}

In this work, we focus our attention on two random flights in
$\mathbb{R}^d$ indicated by $\{\underline{X}_d(t):t\geq 0\}$, with
$d\geq 2$, and $\{\underline{Y}_d(t):t\geq 0\}$, with $d\geq 3$,
respectively. Now we briefly describe these models; for more
details the reader can consult \cite{DO}. The two random motions
start at the origin, move with finite velocity $c$ and, when $n$
changes of direction are recorded, the $n+1$ time lengths of the
displacements performed by the motions have suitable rescaled Dirichlet
distributions. More precisely, let us indicate by $\tau_j,$ with $j\in\{1,...,n+1\}$ and $\tau_0=0$, $\tau_{n+
1}=t-\sum_{j=1}^n\tau_j$, the time displacements of the motions. We assume that
$$f_X(\tau_1,\ldots,\tau_n)=\frac{\Gamma((n+1)(d-1))}{(\Gamma(d-1))^{n+1}}\frac{1}{t^{(n+1)(d-1)-1}}\prod_{j=1}^{n+1}\tau_j^{d-2}$$
for $\{\underline{X}_d(t):t\geq 0\}$ (the rescaled Dirichlet
distributions with parameters $(d-1,\ldots,d-1)$; see eq. (1.4) in
\cite{DO}) and
$$f_Y(\tau_1,\ldots,\tau_n)=\frac{\Gamma((n+1)(\frac{d}{2}-1))}{(\Gamma(\frac{d}{2}-1))^{n+1}}\frac{1}{t^{(n+1)(\frac{d}{2}-1)-1}}
\prod_{j=1}^{n+1}\tau_j^{\frac{d}{2}-2}$$ for
$\{\underline{Y}_d(t):t\geq 0\}$ (the rescaled Dirichlet
distributions with parameters
$\left(\frac{d}{2}-1,\ldots,\frac{d}{2}-1\right)$; see eq. (1.5)
in \cite{DO}) where, in both cases,
$0<\tau_j<t-\sum_{k=0}^{j-1}\tau_k$ for $j\in\{1,\ldots,n\}$. Moreover, in
both cases, the orientations concerning the displacements - the
initial one, and the ones after each change of direction - are
uniformly chosen on the sphere of $\mathbb{R}^d$ with radius one
(see eq. (1.3) in \cite{DO}). Any direction is chosen
independently from the previous one.  The sample paths of
$\{\underline{X}_d(t):t\geq 0\}$, and $\{\underline{Y}_d(t):t\geq
0\}$ appear as random straight-lines with sharp turns. Figure
\ref{fig1} displays a typical sample path.

The aim of this work is to investigate the asymptotic behavior of
the above models as $t\to\infty$. In particular we obtain large
deviation results concerning the random flights. The large
deviation theory provides an asymptotic computation of small
probabilities on exponential scale; estimates based on large
deviations play a crucial role in resolving a variety of questions
in statistics, engineering, statistical mechanics and applied
probability.

The paper is organized as follows. In Section
\ref{sec:preliminaries}, we provide a brief introduction on the
large deviation theory and recall the probability distributions
for the random flights considered in this work. The main results
are contained in Section \ref{sec:results}. In particular we
obtain the large deviation results for the random flights
$\{\underline{X}_d(t):t\geq 0\}$, with $d\geq 2$, and
$\{\underline{Y}_d(t):t\geq 0\}$, with $d\geq 3$. If we consider a
homogeneous Poisson process governing the changes of direction of
the motion, we obtain a standard random flight. We have that
$\{\underline{X}_2(t):t\geq 0\}$ and $\{\underline{Y}_4(t):t\geq
0\}$ represent the standard cases. Section \ref{sec:results}
contains some propositions concerning the standard random flights.
In this last case we are also able to consider the large
deviations results for non-conditional laws. Many comparisons
concerning rate functions of the standard random flights are
listed in Section \ref{sec:final}.

\section{Preliminaries}\label{sec:preliminaries}
We start with some preliminaries on large deviations. We recall
the basic definitions in \cite{DZ} (pages 4--5). Let $\mathcal{Z}$
be a Hausdorff topological space with Borel $\sigma$-algebra
$\mathcal{B}_{\mathcal{Z}}$. A lower semi-continuous function
$I:\mathcal{Z}\to [0,\infty]$ is called rate function. A family of
probability measures $\{\pi_t:t>0\}$ on
$(\mathcal{Z},\mathcal{B}_{\mathcal{Z}})$ satisfies the
\emph{large deviation principle} (LDP for short), as $t\to\infty$,
with rate function $I$ if
$$\limsup_{t\to\infty}\frac{1}{t}\log\pi_t(F)\leq-\inf_{z\in
F}I(z)\quad\textrm{for all closed sets}\ F$$ and
$$\liminf_{t\to\infty}\frac{1}{t}\log\pi_t(G)\geq-\inf_{z\in
G}I(z)\quad\textrm{for all open sets}\ G.$$ A rate function $I$ is
said to be good if all the level sets
$\{\{z\in\mathcal{Z}:I(z)\leq\gamma\}:\gamma\geq 0\}$ are compact.
In what follows we use condition (b) with equation (1.2.8) in
\cite{DZ}, which is equivalent to the lower bound for open sets:
\begin{equation}\label{eq:LB-local-condition}
\liminf_{t\to\infty}\frac{1}{t}\log\pi_t(G)\geq-I(z)\quad
\left.\begin{array}{ll}
\mbox{for all $z\in\mathcal{Z}$ such that $I(z)<\infty$ and}\\
\mbox{for all open sets $G$ such that $z\in G$}.
\end{array}\right.
\end{equation}

\begin{Remark}\label{rem:comparison}
Assume that $\{\pi_t^{(1)}:t>0\}$ and $\{\pi_t^{(2)}:t>0\}$
satisfy the LDP with the rate functions $I_1$ and $I_2$,
respectively; moreover assume that $I_1$ and $I_2$ uniquely vanish
at the same point $z_0$. Then, if we have $I_1(z)>I_2(z)$ for all
$z$ in a neighborhood $U$ of $z_0$ (except $z_0$ because
$I_1(z_0)=I_2(z_0)=0$), we can say that $\{\pi_t^{(1)}:t>0\}$
converges to $z_0$ faster than $\{\pi_t^{(2)}:t>0\}$, as
$t\to\infty$. Indeed, for all $\varepsilon>0$, there exists
$t_\varepsilon$ such that
$\frac{\pi_t^{(1)}(U^c)}{\pi_t^{(2)}(U^c)}\leq
e^{-t(I_1(U^c)-I_2(U^c)+2\varepsilon)}$ for all $t>t_\varepsilon$,
where $I_k(U^c)=\inf_{z\in U^c}I_k(z)$ for $k\in\{1,2\}$; thus,
since $I_1(U^c)>I_2(U^c)>0$, we have
$\frac{\pi_t^{(1)}(U^c)}{\pi_t^{(2)}(U^c)}\to 0$ as $t\to\infty$.
\end{Remark}

Throughout the paper we use the following notation:
$\|\underline{z}_d\|:=\sqrt{z_1^2+\cdots+z_d^2}$ is the norm of
$\underline{z}_d=(z_1,\ldots,z_d)\in\mathbb{R}^d$;
$B_\delta(\underline{z}_d):=\{\underline{y}_d\in\mathbb{R}^d:\|\underline{y}_d-\underline{z}_d\|<\delta\}$
is the neighborhood of $\underline{z}_d$ with radius $\delta>0$;
$\underline{0}_d$ is the null vector in $\mathbb{R}^d$.

Now we recall some preliminaries on the densities obtained in
Theorem 2 in \cite{DO}. They concern two random flights in
$\mathbb{R}^d$: $\{\underline{X}_d(t):t\geq 0\}$ for $d\geq 2$;
$\{\underline{Y}_d(t):t\geq 0\}$ for $d\geq 3$. Moreover the
number $n\geq 1$ of changes of direction in the time interval
$[0,t]$ is fixed and, from now on, we use the following notation:
$$\mu_d(E,t;n):=P_n(\underline{X}_d(t)\in E);\quad \nu_d(E,t;n):=P_n(\underline{Y}_d(t)\in E).$$

Then we have:
$$\mu_d(E,t;n):=\int_{E\cap B_{ct}(\underline{0}_d)}
\frac{\Gamma(\frac{n+1}{2}(d-1)+\frac{1}{2})}{\Gamma(\frac{n}{2}(d-1))}
\frac{(c^2t^2-\|\underline{x}_d\|^2)^{\frac{n}{2}(d-1)-1}}{\pi^{d/2}(ct)^{(n+1)(d-1)-1}}
dx_1\cdots dx_d\quad (d\geq 2);$$
$$\nu_d(E,t;n):=\int_{E\cap B_{ct}(\underline{0}_d)}
\frac{\Gamma((n+1)(\frac{d}{2}-1)+1)}{\Gamma(n(\frac{d}{2}-1))}
\frac{(c^2t^2-\|\underline{y}_d\|^2)^{n(\frac{d}{2}-1)-1}}{\pi^{d/2}(ct)^{2(n+1)(\frac{d}{2}-1)}}
dy_1\cdots dy_d\quad (d\geq 3).$$ Note that in both cases we have
a probability measure $\xi_d(\cdot,t;n)$ defined by
\begin{equation}\label{eq:1-formula-for-2-cases}
\xi_d(E,t;n)=\int_{E\cap
B_{ct}(\underline{0}_d)}h_d(\underline{z}_d,t;n)dz_1\cdots dz_d
\end{equation}
for some density $h_d(\cdot,t;n)$ which has the universal
isotropic form
\begin{equation}\label{eq:1-density-for-2-cases}
h_d(\underline{z}_d,t;n)=\alpha(n)t^{-\gamma(n)}(c^2t^2-\|\underline{z}_d\|^2)^{\beta(n)}\
(\mathrm{for}\ \underline{z}_d\in B_{ct}(\underline{0}_d)),
\end{equation}
more precisely we have
\begin{equation}\label{eq:details-for-each-case}
\left\{\begin{array}{lllll}
\underline{X}_d&\rightsquigarrow
&\gamma(n)=(n+1)(d-1)-1,&\alpha(n)=\frac{\Gamma(\frac{n+1}{2}(d-1)+\frac{1}{2})}{\Gamma(\frac{n}{2}(d-1))\pi^{d/2}c^{\gamma(n)}},
&\beta(n)=\frac{n}{2}(d-1)-1;\\
\underline{Y}_d&\rightsquigarrow
&\gamma(n)=2(n+1)(\frac{d}{2}-1),&\alpha(n)=\frac{\Gamma((n+1)(\frac{d}{2}-1)+1)}{\Gamma(n(\frac{d}{2}-1))\pi^{d/2}c^{\gamma(n)}},
&\beta(n)=n(\frac{d}{2}-1)-1.
\end{array}\right.
\end{equation}

A standard random flight $\{Z_d(t):t\geq 0\}$ in $\mathbb{R}^d$ is
a random motion which starts at the origin $\underline{0}_d$,
moves with constant velocity $c$, chooses the directions uniformly
on the $d$-dimensional sphere with radius 1, and changes direction
at any occurrence of a homogeneous Poisson process $\{N(t):t\geq
0\}$ with intensity $\lambda$. This last assumption implies that
the joint distribution of the lengths between two consecutive
changes of direction is uniform. A reference for this random
motion is \cite{OD}; the case $d=2$ was studied in some earlier
papers (see e.g. section 2 in \cite{S} where $c=1$).

One can check that some laws presented above are suitable
conditional distributions for standard random flights; more
precisely, for each fixed $t>0$, we have
$$\mu_2(\cdot,t;n)=P(Z_2(t)\in\cdot|N(t)=n)\quad \mathrm{and}\quad \nu_4(\cdot,t;n)=P(Z_4(t)\in\cdot|N(t)=n).$$
Indeed, for $d=2$ and $d=4$ the Dirichlet distributions reduce to
the uniform distribution.

Finally, in view of what follows, it is useful to recall the
following limits (which can be proved by inspection).

\begin{Lemma}\label{lemma}
If we have $\lim_{t\to\infty}w_t=w$ for some $w\in(0,\infty)$,
then
$$b(w):=\lim_{t\to\infty}\frac{\beta(tw_t)}{t}=\left\{\begin{array}{ll}
\frac{w}{2}(d-1)&\ if\ \xi_d(\cdot,t;n)=\mu_d(\cdot,t;n)\\
w(\frac{d}{2}-1)&\ if\ \xi_d(\cdot,t;n)=\nu_d(\cdot,t;n)
\end{array}\right.$$
and $\lim_{t\to\infty}\frac{\gamma(tw_t)}{t}=2b(w)$.
\end{Lemma}

\section{Results}\label{sec:results}
We start with Proposition \ref{prop:LDP-conditional} which
provides the LDP for the families of laws
$\{\mu_d(t\cdot,t;n_t):t>0\}$ and $\{\nu_d(t\cdot,t;n_t):t>0\}$,
as $t\to\infty$, under a suitable limit condition on $n_t$.

\begin{Proposition}\label{prop:LDP-conditional}
Assume to have $\lim_{t\to\infty}w_t=w$ for some $w\in(0,\infty)$.
Then:\\
$(a)$ $\{\mu_d(t\cdot,t;tw_t):t>0\}$ satisfies the LDP with good
rate function $I_d(\cdot;w)$ defined by
$$I_d(\underline{z}_d;w)=\left\{\begin{array}{ll}
w(d-1)\log\left(\frac{c}{\sqrt{c^2-\|\underline{z}_d\|^2}}\right)&\ if\ \|\underline{z}_d\|<c\\
\infty&\ otherwise;
\end{array}\right.$$
$(b)$ $\{\nu_d(t\cdot,t;tw_t):t>0\}$ satisfies the LDP with good
rate function $J_d(\cdot;w)$ defined by
$$J_d(\underline{z}_d;w)=\left\{\begin{array}{ll}
2w\left(\frac{d}{2}-1\right)\log\left(\frac{c}{\sqrt{c^2-\|\underline{z}_d\|^2}}\right)&\ if\ \|\underline{z}_d\|<c\\
\infty&\ otherwise.
\end{array}\right.$$
\end{Proposition}
\noindent\emph{Proof.} We consider $\{\xi_d(t\cdot,t;n_t):t>0\}$
defined by
\eqref{eq:1-formula-for-2-cases}-\eqref{eq:1-density-for-2-cases},
and we take into account \eqref{eq:details-for-each-case} for each
case. Thus the rate functions $I_d(\cdot;w)$ and $J_d(\cdot;w)$ in
the statement coincide with
$$K_d(\underline{z}_d;w)=\left\{\begin{array}{ll}
2b(w)\log\left(\frac{c}{\sqrt{c^2-\|\underline{z}_d\|^2}}\right)&\ \mathrm{if}\ \|\underline{z}_d\|<c\\
\infty&\ \mathrm{otherwise},
\end{array}\right.$$
where $b(w)$ is as in Lemma \ref{lemma}. We remark that, by Lemma
\ref{lemma}, we have
$\lim_{t\to\infty}\frac{1}{t}\log\alpha(tw_t)=-2b(w)\log c$. The
proof is divided in two parts.\\
\emph{1) Proof of the lower bound for open sets.} We want to check
the equivalent condition \eqref{eq:LB-local-condition}, i.e.
$$\liminf_{t\to\infty}\frac{1}{t}\log\xi_d(tG,t;tw_t)\geq-K_d(\underline{z}_d;w)$$
for all $\underline{z}_d\in\mathbb{R}^d$ such that
$\|\underline{z}_d\|<c$ and for all open sets $G$ such that
$\underline{z}_d\in G$. Firstly we can take $\varepsilon>0$ small
enough to have $B_\varepsilon(\underline{z}_d)\subset G\cap
B_c(\underline{0}_d)$. Then we have
\begin{align*}
\xi_d(tG,t;tw_t)\geq&\xi_d(tB_\varepsilon(\underline{z}_d),t;tw_t)=\xi_d(B_{\varepsilon
t}(\underline{z}_dt),t;tw_t)\\
=&\int_{B_{\varepsilon
t}(\underline{z}_dt)}\alpha(tw_t)t^{-\gamma(tw_t)}(c^2t^2-\|\underline{v}_d\|^2)^{\beta(tw_t)}dv_1\cdots
dv_d\\
=&\int_{B_\varepsilon(\underline{z}_d)}\alpha(tw_t)t^{-\gamma(tw_t)}(c^2t^2-\|\underline{v}_d\|^2t^2)^{\beta(tw_t)}t^ddv_1\cdots
dv_d\\
\geq&\alpha(tw_t)t^{-\gamma(tw_t)+2\beta(tw_t)+d}(c^2-\sup\{\|\underline{v}_d\|^2:\underline{v}_d\in
B_\varepsilon(\underline{z}_d)\})^{\beta(tw_t)}\mathrm{measure}(B_\varepsilon(\underline{z}_d)),
\end{align*}
whence we obtain
$$\liminf_{t\to\infty}\frac{1}{t}\log\xi_d(tG,t;tw_t)\geq-2b(w)\log
c+b(w)\log(c^2-\sup\{\|\underline{v}_d\|^2:\underline{v}_d\in
B_\varepsilon(\underline{z}_d)\}),$$ and we conclude
by letting $\varepsilon$ go to zero.\\
\emph{2) Proof of the upper bound for closed sets.} We have to
check
$$\limsup_{t\to\infty}\frac{1}{t}\log\xi_d(tF,t;tw_t)\leq-\inf_{\underline{z}_d\in
F}K_d(\underline{z}_d;w)\quad\textrm{for all closed sets}\ F.$$
Firstly note that this condition trivially holds if
$\underline{0}_d\in F$ and if $F\cap
B_c(\underline{0}_d)=\emptyset$. Thus, from now on, we assume that
$\underline{0}_d\notin F$ and $F\cap
B_c(\underline{0}_d)\neq\emptyset$. We can find
$\underline{z}_d^F\in F$ such that
$\|\underline{z}_d^F\|=\inf\{\|\underline{z}_d\|:\underline{z}_d\in
F\cap B_c(\underline{0}_d)\}$; note that
$r_F:=\|\underline{z}_d^F\|\in (0,c)$. Then, since $F\subset
(B_{r_F}(\underline{0}_d))^c$, we have
\begin{align*}
\xi_d(tF,t;tw_t)&\leq\xi_d((B_{r_Ft}(\underline{0}_d))^c,t;tw_t)
=\int_{r_Ft}^{ct}\alpha(tw_t)t^{-\gamma(tw_t)}(c^2t^2-\rho^2)^{\beta(tw_t)}\frac{2\pi^{d/2}}{\Gamma(d/2)}\rho^{d-1}d\rho\\
&\leq\alpha(tw_t)t^{-\gamma(tw_t)}\frac{\pi^{d/2}}{\Gamma(d/2)}(ct)^{d-2}\int_{r_Ft}^{ct}(c^2t^2-\rho^2)^{\beta(tw_t)}2\rho d\rho\\
&=\alpha(tw_t)t^{-\gamma(tw_t)}\frac{\pi^{d/2}}{\Gamma(d/2)}(ct)^{d-2}
\left[-\frac{(c^2t^2-\rho^2)^{\beta(tw_t)+1}}{\beta(tw_t)+1}\right]_{\rho=r_Ft}^{\rho=ct}\\
&=\alpha(tw_t)t^{-\gamma(tw_t)+2\beta(tw_t)+2}\frac{\pi^{d/2}}{\Gamma(d/2)}(ct)^{d-2}
\frac{(c^2-r_F^2)^{\beta(tw_t)+1}}{\beta(tw_t)+1}
\end{align*}
whence we obtain
\begin{align*}
\limsup_{t\to\infty}\frac{1}{t}\log\xi_d(tF,t;tw_t)&\leq-2b(w)\log
c+b(w)\log(c^2-r_F^2)=-2b(w)\log\left(\frac{c}{\sqrt{c^2-r_F^2}}\right)\\
&=-2b(w)\log\left(\frac{c}{\sqrt{c^2-\|\underline{z}_d^F\|^2}}\right)
=-K_d(\underline{z}_d^F;w).
\end{align*}
In conclusion, since
$K_d(\underline{z}_d;w)\geq(\mathrm{resp.}=)K_d(\underline{u}_d;w)$
if and only if
$\|\underline{z}_d\|\geq(\mathrm{resp.}=)\|\underline{u}_d\|$, we
have $K_d(\underline{z}_d^F;w)=\inf_{\underline{z}_d\in
F}K_d(\underline{z}_d;w)$ and this completes the proof. $\Box$\\

The next Propositions
\ref{prop:LDP-planar-random-motion}-\ref{prop:LDP-R4-random-motion}
provide the LDP for the non-conditional laws
$\{P(Z_2(t)\in\cdot):t>0\}$ and $\{P(Z_4(t)\in\cdot):t>0\}$ and,
in analogy with the symbols introduced above, from now on we use
the following notation:
$$\mu_2(\cdot,t)=P(Z_2(t)\in\cdot);\quad\nu_4(\cdot,t)=P(Z_4(t)\in\cdot).$$

We start with the first family of laws: for each fixed $t>0$, it
is known (see e.g. eq. (1.2) in \cite{OD}) that $\mu_2(\cdot,t)$
has an absolutely continuous part given by
\begin{equation}\label{eq:ac-part-d=2}
\frac{\lambda}{2\pi c}\frac{e^{-\lambda
t+\frac{\lambda}{c}\sqrt{c^2t^2-\|\underline{z}_2\|^2}}}{\sqrt{c^2t^2-\|\underline{z}_2\|^2}}
1_{B_{ct}(\underline{0}_2)}(\underline{z}_2)dz_1dz_2,
\end{equation}
and a singular part uniformly distributed on the boundary of
$B_{ct}(\underline{0}_2)$ with weight $e^{-\lambda t}$.

\begin{Proposition}\label{prop:LDP-planar-random-motion}
The family $\{\mu_2(t\cdot,t):t>0\}$ satisfies the LDP with good
rate function $I_2$ defined by
$$I_2(\underline{z}_2)=\left\{\begin{array}{ll}
\lambda\left(1-\sqrt{1-\frac{\|\underline{z}_2\|^2}{c^2}}\right)&\ if\ \|\underline{z}_2\|\leq c\\
\infty&\ otherwise.
\end{array}\right.$$
\end{Proposition}
\noindent\emph{Proof.} The proof is divided in two parts.\\
\emph{1) Proof of the lower bound for open sets.} We want to check
the equivalent condition \eqref{eq:LB-local-condition}, i.e.
$$\liminf_{t\to\infty}\frac{1}{t}\log\mu_2(tG,t)
\geq-\lambda\left(1-\sqrt{1-\frac{\|\underline{z}_2\|^2}{c^2}}\right)$$
for all $\underline{z}_2\in\mathbb{R}^2$ such that
$\|\underline{z}_2\|\leq c$ and for all open sets $G$ such that
$\underline{z}_2\in G$. Firstly we can take $\varepsilon>0$ small
enough to have $B_\varepsilon(\underline{z}_2)\subset G$;
moreover, if $\underline{z}_2\in B_c(\underline{0}_2)$, we also
require that $B_\varepsilon(\underline{z}_2)\subset
B_c(\underline{0}_2)$. Then we have
\begin{align*}
\mu_2(tG,t)\geq&\mu_2(tB_\varepsilon(\underline{z}_2),t)
=\mu_2(B_{\varepsilon t}(\underline{z}_2t),t)\\
\geq&\int_{B_{\varepsilon t}(\underline{z}_2t)\cap
B_{ct}(\underline{0}_2)}\frac{\lambda}{2\pi c}\frac{e^{-\lambda
t+\frac{\lambda}{c}\sqrt{c^2t^2-\|\underline{v}_2\|^2}}}{\sqrt{c^2t^2-\|\underline{v}_2\|^2}}dv_1dv_2\\
\geq&\int_{B_\varepsilon(\underline{z}_2)\cap
B_c(\underline{0}_2)}\frac{\lambda}{2\pi c}\frac{e^{-\lambda
t+\frac{\lambda t}{c}\sqrt{c^2-\|\underline{v}_2\|^2}}}{ct}t^2dv_1dv_2\\
\geq&\frac{\lambda t}{2\pi c}\frac{e^{-\lambda t+\lambda
t\sqrt{1-\frac{\sup\{\|\underline{v}_2\|^2:\underline{v}_2\in
B_\varepsilon(\underline{z}_2)\cap
B_c(\underline{0}_2)\}}{c^2}}}}{c}\cdot
\underbrace{\mathrm{measure}(B_\varepsilon(\underline{z}_2)\cap
B_c(\underline{0}_2))}_{>0},
\end{align*}
and, arguing as in the proof of Proposition
\ref{prop:LDP-conditional}, we conclude by taking
$\liminf_{t\to\infty}\frac{1}{t}\log$ (for both the left hand side
and the right hand side) and by letting $\varepsilon$ go to
zero.\\
\emph{2) Proof of the upper bound for closed sets.} We have to
check
$$\limsup_{t\to\infty}\frac{1}{t}\log
\mu_2(tF,t)\leq-\inf_{\underline{z}_2\in
F}I_2(\underline{z}_2)\quad\textrm{for all closed sets}\ F.$$
Firstly note that this condition trivially holds if
$\underline{0}_2\in F$ and if
$F\cap\overline{B_c(\underline{0}_2)}=\emptyset$. Thus, from now
on, we assume that $\underline{0}_2\notin F$ and
$F\cap\overline{B_c(\underline{0}_2)}\neq\emptyset$. We can find
$\underline{z}_2^F\in F$ such that
$\|\underline{z}_2^F\|=\inf\{\|\underline{z}_2\|:\underline{z}_2\in
F\cap\overline{B_c(\underline{0}_2)}\}$; note that
$r_F:=\|\underline{z}_2^F\|\in (0,c]$. Then, since $F\subset
(B_{r_F}(\underline{0}_2))^c$, we have
\begin{align*}
\mu_2(tF,t)&\leq\mu_2((B_{r_Ft}(\underline{0}_2))^c,t)=\int_{r_Ft}^{ct}\frac{\lambda}{2\pi
c}\frac{e^{-\lambda
t+\frac{\lambda}{c}\sqrt{c^2t^2-\rho^2}}}{\sqrt{c^2t^2-\rho^2}}2\pi\rho
d\rho+e^{-\lambda t}\\
&=\left[-e^{-\lambda
t+\frac{\lambda}{c}\sqrt{c^2t^2-\rho^2}}\right]_{\rho=r_Ft}^{\rho=ct}+e^{-\lambda
t}=\exp\left(-\lambda
t+\frac{\lambda}{c}\sqrt{c^2t^2-r_F^2t^2}\right),
\end{align*}
whence we obtain
$$\limsup_{t\to\infty}\frac{1}{t}\log
\mu_2(tF,t)\leq-\lambda+\frac{\lambda}{c}\sqrt{c^2-r_F^2}=
-\lambda\left(1-\sqrt{1-\frac{\|\underline{z}_2^F\|^2}{c^2}}\right)=-I_2(\underline{z}_2^F).$$
Thus, arguing as in the proof of Proposition
\ref{prop:LDP-conditional}, we conclude noting that
$I_2(\underline{z}_2^F)=\inf_{\underline{z}_2\in
F}I_2(\underline{z}_2)$ because
$I_2(\underline{z}_2)\geq(\mathrm{resp.}=)I_2(\underline{u}_2)$ if
and only if
$\|\underline{z}_2\|\geq(\mathrm{resp.}=)\|\underline{u}_2\|$.
$\Box$\\

Now we consider the second family of laws: for each fixed $t>0$,
it is known (see e.g. Theorem 3.2 in \cite{OD}) that
$\nu_4(\cdot,t)$ has an absolutely continuous part given by
\begin{equation}\label{eq:ac-part-d=4}
\frac{\lambda}{c^4t^3\pi^2}e^{-\frac{\lambda}{c^2t}\|\underline{z}_4\|^2}
\left(2+\frac{\lambda}{c^2t}\left(c^2t^2-\|\underline{z}_4\|^2\right)\right)
1_{B_{ct}(\underline{0}_4)}(\underline{z}_4)dz_1dz_2dz_3dz_4,
\end{equation}
and a singular part uniformly distributed on the boundary of
$B_{ct}(\underline{0}_4)$ with weight $e^{-\lambda t}$.

\begin{Proposition}\label{prop:LDP-R4-random-motion}
The family $\{\nu_4(t\cdot,t):t>0\}$ satisfies the LDP with good
rate function $J_4$ defined by
$$J_4(\underline{z}_4)=\left\{\begin{array}{ll}
\frac{\lambda}{c^2}\|\underline{z}_4\|^2&\ if\ \|\underline{z}_4\|\leq c\\
\infty&\ otherwise.
\end{array}\right.$$
\end{Proposition}
\noindent\emph{Proof.} The proof is divided in two parts.\\
\emph{1) Proof of the lower bound for open sets.} We want to check
the equivalent condition \eqref{eq:LB-local-condition}, i.e.
$$\liminf_{t\to\infty}\frac{1}{t}\log\nu_4(tG,t)\geq-\frac{\lambda}{c^2}\|\underline{z}_4\|^2$$
for all $\underline{z}_4\in\mathbb{R}^4$ such that
$\|\underline{z}_4\|\leq c$ and for all open sets $G$ such that
$\underline{z}_4\in G$. Firstly we can take $\varepsilon>0$ small
enough to have $B_\varepsilon(\underline{z}_4)\subset G$;
moreover, if $\underline{z}_4\in B_c(\underline{0}_4)$, we also
require that $B_\varepsilon(\underline{z}_4)\subset
B_c(\underline{0}_4)$. Then we have
\begin{align*}
\nu_4(tG,t)\geq&\nu_4(tB_\varepsilon(\underline{z}_4),t)
=\nu_4(B_{\varepsilon t}(\underline{z}_4t),t)\\
\geq&\int_{B_{\varepsilon t}(\underline{z}_4t)\cap
B_{ct}(\underline{0}_4)}\frac{\lambda}{c^4t^3\pi^2}e^{-\frac{\lambda}{c^2t}\|\underline{v}_4\|^2}
\left(2+\frac{\lambda}{c^2t}\left(c^2t^2-\|\underline{v}_4\|^2\right)\right)dv_1dv_2dv_3dv_4\\
=&\int_{B_\varepsilon(\underline{z}_4)\cap
B_c(\underline{0}_4)}\frac{\lambda}{c^4t^3\pi^2}e^{-\frac{\lambda
t}{c^2}\|\underline{v}_4\|^2}
\left(2+\frac{\lambda t}{c^2}\left(c^2-\|\underline{v}_4\|^2\right)\right)t^4dv_1dv_2dv_3dv_4\\
\geq&\frac{\lambda t}{c^4\pi^2}e^{-\frac{\lambda
t}{c^2}\sup\{\|\underline{v}_4\|^2:\underline{v}_4\in
B_\varepsilon(\underline{z}_4)\cap B_c(\underline{0}_4)\}}\\
&\cdot\left(2+\frac{\lambda
t}{c^2}\left(c^2-\sup\{\|\underline{v}_4\|^2:\underline{v}_4\in
B_\varepsilon(\underline{z}_4)\cap
B_c(\underline{0}_4)\}\right)\right)\underbrace{\mathrm{measure}(B_\varepsilon(\underline{z}_4)\cap
B_c(\underline{0}_4))}_{>0},
\end{align*}
and we conclude following the lines of Proposition
\ref{prop:LDP-planar-random-motion} (final part of the proof of
the lower bound).\\
\emph{2) Proof of the upper bound for closed sets.} We have to
check
$$\limsup_{t\to\infty}\frac{1}{t}\log
\nu_4(tF,t)\leq-\inf_{\underline{z}_4\in
F}J_4(\underline{z}_4)\quad\textrm{for all closed sets}\ F.$$
Firstly note that this condition trivially holds if
$\underline{0}_4\in F$ and if
$F\cap\overline{B_c(\underline{0}_4)}=\emptyset$. Thus, from now
on, we assume that $\underline{0}_4\notin F$ and
$F\cap\overline{B_c(\underline{0}_4)}\neq\emptyset$. We can find
$\underline{z}_4^F\in F$ such that
$\|\underline{z}_4^F\|=\inf\{\|\underline{z}_4\|:\underline{z}_4\in
F\cap\overline{B_c(\underline{0}_4)}\}$; note that
$r_F:=\|\underline{z}_4^F\|\in (0,c]$. Then, since $F\subset
(B_{r_F}(\underline{0}_4))^c$, we have
$$\nu_4(tF,t)\leq\nu_4((B_{r_Ft}(\underline{0}_4))^c,t)=\int_{r_Ft}^{ct}\frac{\lambda}{c^4t^3\pi^2}e^{-\frac{\lambda}{c^2t}\rho^2}
\left(2+\frac{\lambda}{c^2t}(c^2t^2-\rho^2)\right)2\pi^2\rho^3d\rho+e^{-\lambda
t};$$
moreover we note that, if $\rho\in[r_Ft,ct]$, we have
$$\left\{\begin{array}{ll}
0\leq\frac{\lambda}{c^2t}(c^2t^2-\rho^2)=\lambda
t\left(1-\left(\frac{\rho}{ct}\right)^2\right)\leq\lambda t\\
0\leq\frac{2\lambda\rho^3}{c^4t^3}=2\left(\frac{\rho}{ct}\right)^2\frac{\lambda\rho}{c^2t}\leq\frac{2\lambda\rho}{c^2t},
\end{array}\right.$$
whence we obtain
\begin{align*}
\nu_4(tF,t)\leq&(2+\lambda
t)\int_{r_Ft}^{ct}\frac{2\lambda\rho}{c^2t}
e^{-\frac{\lambda}{c^2t}\rho^2}d\rho+e^{-\lambda t}\\
=&(2+\lambda
t)\left[-e^{-\frac{\lambda}{c^2t}\rho^2}\right]_{\rho=r_Ft}^{\rho=ct}+e^{-\lambda
t}
=(2+\lambda t)e^{-\lambda\frac{r_F^2}{c^2}t}\left(1-e^{-\lambda
t\left(1-\frac{r_F^2}{c^2}\right)}\right)+e^{-\lambda t};
\end{align*}
actually, if $r_F=c$, we have
$\nu_4((B_{r_Ft}(\underline{0}_4))^c,t)=P(N(t)=0)=e^{-\lambda t}$.
In conclusion, by Lemma 1.2.15 in \cite{DZ} (we need this for
$r_F\in(0,c)$; the case $r_F=c$ is trivial), we obtain
$$\limsup_{t\to\infty}\frac{1}{t}\log
\nu_4(tF,t)\leq
-\frac{\lambda}{c^2}r_F^2=-J_4(\underline{z}_4^F);$$ then we can
conclude following the lines of Proposition
\ref{prop:LDP-planar-random-motion} (final part of the proof of
the upper bound). $\Box$\\

We remark that the rate function $J_4$ in the previous Proposition
\ref{prop:LDP-R4-random-motion} is quadratic on
$\overline{B_c(\underline{0}_4)}$; quadratic rate functions
typically come up in the case of LDPs for Gaussian random
variables and here, if we neglect the factor
$2+\frac{\lambda}{c^2t}\left(c^2t^2-\|\underline{z}_4\|^2\right)$,
the density in \eqref{eq:ac-part-d=4} is proportional to a
Gaussian density restricted on $B_{ct}(\underline{0}_4)$.

\section{Minor results and concluding remarks}\label{sec:final}
We start with the comparison between the rate functions for the
conditional laws and the non-conditional laws related to the
standard random flights $\{Z_2(t):t\geq 0\}$ and $\{Z_4(t):t\geq
0\}$, in the spirit of Remark \ref{rem:comparison}. In both cases
the convergence of the conditional distributions is faster than
the convergence of the non-conditional distributions if and only
if $w\geq\lambda$. For the case $d=2$ we can refer to a result for
the telegraph process on the real line, i.e. Proposition 2.3 in
\cite{DM} (and its consequences), specified to the case
$c_1=c_2=c$: more precisely we have to consider the rate functions
$I_{\lambda,\lambda,c,c}^X$ (for the non-conditional laws) and the
$I_{\lambda,\lambda,c,c}^{X|N}(\cdot;w)$ (for the conditional
laws) in \cite{DM}, and the equalities
\begin{equation}\label{eq:d=2-and-telegraph}
\left\{\begin{array}{l}
I_2(\underline{z}_2)=I_{\lambda,\lambda,c,c}^X(\|\underline{z}_2\|)\\
I_2(\underline{z}_2;w)=I_{\lambda,\lambda,c,c}^{X|N}(\|\underline{z}_2\|;w)\
(\mbox{for all}\ w>0)
\end{array}\right.\ \mbox{for all}\ \underline{z}_2\in\mathbb{R}^2.
\end{equation}
For the case $d=4$ we have the following similar
result.

\begin{figure}[ht]
\begin{center}
\includegraphics[angle=0,width=0.90\textwidth]{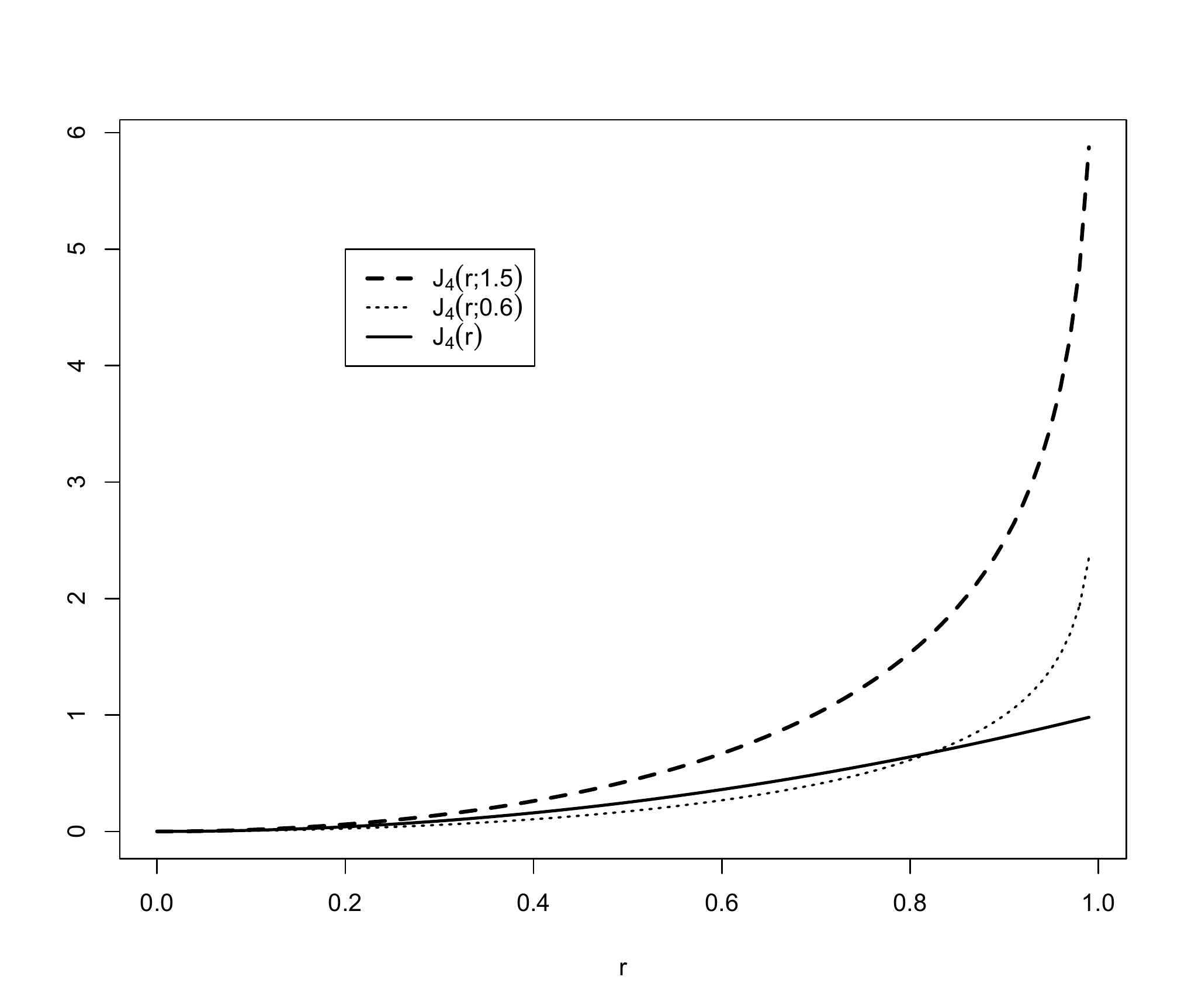}
\caption{The functions $J_4(\underline{z}_4)$ and
$J_4(\underline{z}_4;w)$ for $\lambda=c=1$, where
$r=\|\underline{z}_4\|\in[0,1)$. Two choices for $w$: $w=1.5\geq
1$; $w=0.6\in(0,1)$.}\label{fig2}
\end{center}
\end{figure}

\begin{Proposition}\label{prop:comparison-d=4}
We have two cases. (i) For $w\geq\lambda$, we have
$J_4(\underline{z}_4)\leq J_4(\underline{z}_4;w)$ for all
$\underline{z}_4\in\mathbb{R}^4$; moreover the inequality is
strict for $\underline{z}_4\in
\overline{B_c(\underline{0}_4)}\setminus\{\underline{0}_4\}$. (ii)
For $w\in(0,\lambda)$, there exists $\gamma\in (\xi,1)$ where
$\xi=\sqrt{1-\frac{w}{\lambda}}$ such that:
$J_4(\underline{z}_4;w)>J_4(\underline{z}_4)$ for
$\|\underline{z}_4\|\in(\gamma c,c]$,
$J_4(\underline{z}_4;w)<J_4(\underline{z}_4)$ for
$\|\underline{z}_4\|\in(0,\gamma c)$ and
$J_4(\underline{z}_4;w)=J_4(\underline{z}_4)$ otherwise.
\end{Proposition}
\noindent Proposition \ref{prop:comparison-d=4} can be proved by
inspection and we omit the details. Some inequalities in
Proposition \ref{prop:comparison-d=4} between $J_4$ and
$J_4(\cdot;w)$ are displayed in Figure \ref{fig2} below when
$\lambda=c=1$.

It is known (see the discussion in \cite{OD}) that the telegraph
process without drift on the line has a strict analogy with the
standard random flight in $\mathbb{R}^2$ and not in
$\mathbb{R}^4$. Actually the density in \eqref{eq:ac-part-d=2}
satisfies the two-dimensional telegraph equation (see e.g. the
displayed equation just after (1.2) in \cite{OD}). We have a
similar situation when we compare the rate functions for the
telegraph process (without drift) on the line and for the standard
random flights. Firstly we have \eqref{eq:d=2-and-telegraph} for
the case $d=2$. On the contrary, for $d=4$, for all
$\underline{z}_4\in\mathbb{R}^4$ we have:
$$J_4(\underline{z}_4)=\left\{\begin{array}{ll}
\frac{\lambda}{c^2}\|\underline{z}_4\|^2&\ \mathrm{if}\ \|\underline{z}_4\|\leq c\\
\infty&\ \mathrm{otherwise};
\end{array}\right.\quad I_{\lambda,\lambda,c,c}^X(\|\underline{z}_4\|)=\left\{\begin{array}{ll}
\lambda\left(1-\sqrt{1-\frac{\|\underline{z}_4\|^2}{c^2}}\right)&\ \mathrm{if}\ \|\underline{z}_4\|\leq c\\
\infty&\ \mathrm{otherwise};
\end{array}\right.$$
$$J_4(\underline{z}_4;w)=2I_{\lambda,\lambda,c,c}^{X|N}(\|\underline{z}_4\|;w)\ \mbox{for all}\ w>0.$$
Therefore $J_4(\cdot)$ and $I_{\lambda,\lambda,c,c}^X(\|\cdot\|)$
are quite different, while $J_4(\cdot;w)$ and
$I_{\lambda,\lambda,c,c}^{X|N}(\|\cdot\|;w)$ differ for the
multiplicative factor 2.

In view of what follows we remark that the values assumed by the
functions $I_2,I_2(\cdot;w),J_4,J_4(\cdot;w)$ depend on the
distance of their arguments from the origin; therefore we use the
symbol $r$ in place of $\|\underline{z}_2\|$ and
$\|\underline{z}_4\|$, and we write
$I_2(r),I_2(r;w),J_4(r),J_4(r;w)$ instead of
$I_2(\underline{z}_2),I_2(\underline{z}_2;w),J_4(\underline{z}_4),J_4(\underline{z}_4;w)$
with a slight abuse of notation. Then, for all $r\geq 0$, we have
$$\left\{\begin{array}{l}
J_4(r;w)=2I_2(r;w)\geq I_2(r;w)\\
J_4(r)\geq I_2(r)
\end{array}\right.$$
(the first inequality is trivial, the second one can be checked
with easy computations), and Figure \ref{fig3} below displays the
second inequality when $\lambda=c=1$. We also remark that the two
inequalities above turn into an inequality if and only if we have
one the following cases:
$$J_4(r;w)=I_2(r;w)=\left\{\begin{array}{ll}
0&\ \mathrm{if}\ r=0\\
\infty&\ \mathrm{if}\ r\in[c,\infty);
\end{array}\right.\quad
J_4(r)=I_2(r)=\left\{\begin{array}{ll}
0&\ \mathrm{if}\ r=0\\
\lambda&\ \mathrm{if}\ r=c\\
\infty&\ \mathrm{if}\ r\in (c,\infty).
\end{array}\right.$$
In particular we remark that $J_4(c)=I_2(c)=\lambda$ for the
non-conditional laws concerns the case without changes of
direction. In conclusion, by taking into account Remark
\ref{rem:comparison}, the convergence at the origin in
$\mathbb{R}^4$ is faster than the analogous convergence at the
origin in $\mathbb{R}^2$ (for both conditional and non-conditional
laws); in some sense this is not surprising because one can expect
to have a faster convergence of a normalized random flight in a
higher space.

\begin{figure}[ht]
\begin{center}
\includegraphics[angle=0,width=0.90\textwidth]{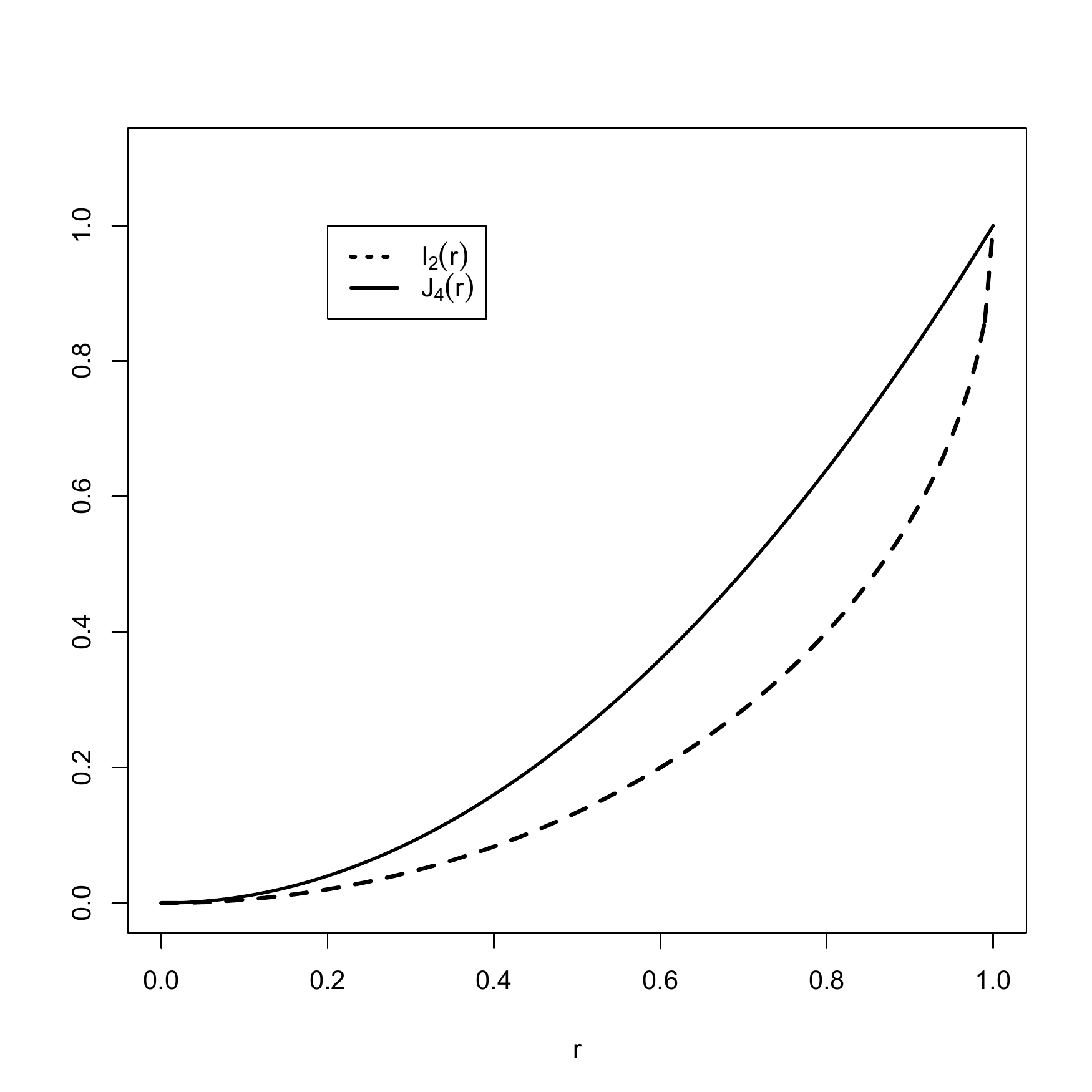}
\caption{The functions $I_2(r)$ and $J_4(r)$ for $\lambda=c=1$,
where $r=\|\underline{z}_d\|\in[0,1]$.}\label{fig3}
\end{center}
\end{figure}

We point out another difference between the cases $d=2$ and $d=4$.
We recall that, if we consider the telegraph process without
drift, i.e. $\lambda_1=\lambda_2=\lambda$ and $c_1=c_2=c$, the
telegraph equation converges to the heat equation as
$\lambda\to\infty$ and $\frac{c^2}{\lambda}\to\sigma^2$ (this and
other connections with the Brownian motion can be found in
\cite{O-1990}, section 4). The corresponding convergence of the
large deviation rate functions under the same scaling is
illustrated in \cite{M-CSSC} (subsection 4.2); more precisely the
convergence of the large deviation rate functions proved there
concerns the more general case with drift, and there are also
convergence results for the decay rates of suitable level crossing
probabilities. Here we remark that, if we consider the same limits
for the rate functions in Propositions
\ref{prop:LDP-planar-random-motion}-\ref{prop:LDP-R4-random-motion}
(again as $\lambda\to\infty$ and
$\frac{c^2}{\lambda}\to\sigma^2$), we have
$I_2(\underline{z}_2)\to\frac{\|\underline{z}_2\|^2}{2\sigma^2}$
and
$J_4(\underline{z}_4)\to\frac{\|\underline{z}_4\|^2}{\sigma^2}$;
then, if we consider the rate function $H_d$ defined by
$H_d(\underline{z}_d)=\frac{\|\underline{z}_d\|^2}{2\sigma^2}$ for
the LDP of $\left\{\frac{B_d(\sigma^2)}{\sqrt{t}}:t>0\right\}$,
where $B_d$ is a standard $d$-dimensional centered Brownian
motion, $H_2(\underline{z}_2)$ coincides with the limit for
$I_2(\underline{z}_2)$, while $H_4(\underline{z}_4)$ is slightly
different from the limit for $J_4(\underline{z}_4)$.

Finally, motivated by potential applications in biology, it is
interesting to present the asymptotic lower bounds in
\eqref{eq:asymptoticLB} for the exit probabilities
$$\Psi_{Z_d}(t;r)=P(\{\|Z_d(s)\|>r\ \mbox{for some}\ s\in[0,t]\})\ \mbox{for}\ d\in\{2,4\}.$$

Actually, for all $r\in(0,c)$, we have
\begin{equation}\label{eq:asymptoticLB}
\liminf_{t\to\infty}\frac{1}{t}\log\Psi_{Z_d}(t;rt)\geq-\ell(d),\
\mbox{where}\ \ell(d):= \left\{\begin{array}{ll}
\lambda\left(1-\sqrt{1-\frac{r^2}{c^2}}\right)&\ \mbox{for}\ d=2\\
\lambda\frac{r^2}{c^2}&\ \mbox{for}\ d=4
\end{array}\right.
\end{equation}
by the inequality $\Psi_{Z_d}(t;rt)\geq
P\left(\frac{\|Z_d(t)\|}{t}>r\right)$ (for all $t,r>0$) and by the
LDPs in Propositions
\ref{prop:LDP-planar-random-motion}-\ref{prop:LDP-R4-random-motion}.

For instance, $r$ could represent a critical threshold in the analysis of the behavior of the bacteria. In some context is realistic to assume that the motile cells moving far from the starting point lose their intrinsic properties. The above lower bound permits us to provide some information about this event.

Obviously, if it is possible to observe the changes of direction
of the random motion, then it would be good to consider the
analogous asymptotic lower bounds by referring to the conditional
laws (thus by considering Proposition \ref{prop:LDP-conditional}
instead of Propositions
\ref{prop:LDP-planar-random-motion}-\ref{prop:LDP-R4-random-motion}).


\end{document}